# Type Two Cuts, Bad Cuts and Very Bad Cuts[1]

Renling Jin[2]

Abstract

Type two cuts, bad cuts and very bad cuts are introduced in [KL] for studying the relationship between Loeb measure and U-topology of a hyperfinite time line in an $\omega_1$-saturated nonstandard universe. The questions concerning the existence of those cuts are asked there. In this paper we answer, fully or partially, some of those questions by showing that: (1) type two cuts exist, (2) the $\aleph_1$-isomorphism property implies that bad cuts exist, but no bad cuts are very bad.

# 0  Introduction

Although the questions in [KL] were discussed in a nonstandard universe, the definition of type two cuts and the proof of their existence in §1 can be carried out in any nonstandard models of Peano Arithmetic. We choose to work within some nonstandard universe throughout the whole paper only for coherence. The reader who is not familiar with nonstandard universe and is only interested in the result in §1 could consider a nonstandard universe as a nonstandard model of Peano Arithmetic in both §0 and §1 without any difficulties. All nonstandard universes mentioned in this paper are $\omega_1$-saturated. Given a nonstandard universe $V$, let $^*\mathbb{N}$ denote the set of all positive integers in $V$ and $\mathbb{N}$ denote the set of all standard positive integers. Following [KL], a non-empty initial segment $U$ of $^*\mathbb{N}$ (under the natural order of $^*\mathbb{N}$) is called a cut if $U$ is closed under addition, i.e. $(\forall x, y \in U)(x + y \in U)$ is true ($U$ is called an additive cut in some literature). For example, $\mathbb{N}$ is the smallest cut and $^*\mathbb{N}$ is the largest cut. There are several ways of constructing new cuts from given cuts shown in [KL]. For example, if $U$ is a cut and $x$ is an element in $^*\mathbb{N}$, then the set

$$xU = \{y \in {}^*\mathbb{N} : (\exists z \in U)(y < xz)\}$$

is a cut. If the element $x$ is in $^*\mathbb{N} \smallsetminus U$, then the set

$$x/U = \{y \in {}^*\mathbb{N} : (\forall z \in U)(y < x/z)\}$$

---

[1]*Mathematics Subject Classification* Primary 03H05, 03H15, 03C50
[2]This research was partially supported by NSF postdoctoral fellowship #DMS-9508887.

is also a cut. Given a cut $U$, let

$$M(U) = \{y \in {}^*\mathbb{N} : (\forall z \in U)\, (yz \in U)\},$$

where $M$ suggests multiplication. Then $M(U)$ is a cut and closed under multiplication. The cut $U$ is called a type one cut if $U = xM(U)$ for some $x \in U$ or $U = x/M(U)$ for some $x \in {}^*\mathbb{N} \smallsetminus U$. $U$ is called a type two cut if it is not type one. As mentioned in [KL], type one–type two cuts are defined in [G]. In §1 we show that type two cuts exist, which answers a question in [KL].

Cuts are used in [KL] for defining $U$-topologies on a hyperfinite time line. Given a hyperinteger $H \in {}^*\mathbb{N} \smallsetminus \mathbb{N}$. The set $\mathcal{H} = \{1, 2, \ldots, H\} \subseteq {}^*\mathbb{N}$ is called a hyperfinite time line. Let $U \subseteq \mathcal{H}$ be a cut. A set $O \subseteq \mathcal{H}$ is called $U$-open if

$$(\forall x \in O)\, (\exists y \in \mathcal{H} \smallsetminus U)\, (\{z \in \mathcal{H} : x - y < z < x + y\} \subseteq O)$$

is true. The $U$-topology of $\mathcal{H}$ is the topology of all $U$-open sets in $\mathcal{H}$. A $U$-topology could also be viewed as an analogue of order topology (note that the natural order topology of $\mathcal{H}$ is discrete). Given $x, y \in \mathcal{H}$, define $a \sim_U y$ iff $|x - y| \in U$. Then it is easy to see that $\sim_U$ is an equivalence relation (here we use the fact that $U$ is closed under addition). Let $x \in \mathcal{H}$. A $\sim_U$ equivalence class containing $x$ is called a $U$-monad of $x$. Given $x, y \in \mathcal{H}$. Define $x \ll_U y$ iff $x < y$ and $x \not\sim_U y$. For any $x, y \in \mathcal{H}$ let

$$I(x, y) = \{z \in \mathcal{H} : x \ll_U z \ll_U y\}.$$

Then the $U$-topology of $\mathcal{H}$ is actually the topology generated by open "intervals" $I(x, y)$ for all $x, y \in \mathcal{H}$. So a $U$-topology is like an order topology by ordering all $U$-monads.

Given a hyperfinite time line $\mathcal{H}$, there is a natural way to define a probability measure called Loeb measure on $\mathcal{H}$. For any internal subset $A$ of $\mathcal{H}$ let $\mu(A) = |A|/H$, where $H$ is the largest number in $\mathcal{H}$. Then $\mu$ is a finite additive, internal uniform counting measure on the algebra of all internal subsets of $\mathcal{H}$. The Loeb measure $L(\mu)$ is now the extension of $st \circ \mu$ to the completion of the $\sigma$-algebra generated by all internal subsets of $\mathcal{H}$, where $st$ is the standard part map. Loeb measure behaves very much like Lebesgue measure on the unit interval $[0, 1]$ of the standard real line.

In [KL] Keisler and Leth probe the similarities between a hyperfinite time line $\mathcal{H}$ equipped with Loeb measure and a $U$-topology, and the standard unit interval $[0, 1]$ equipped with Lebesgue measure and the natural order topology. They consider a cut $U$ to behave nicely if it makes $\mathcal{H}$ much like $[0,1]$. For example, considering the fact that $[0,1]$ contains a meager set of Lebesgue measure one, they call a cut $U \subseteq \mathcal{H}$ a good cut if $\mathcal{H}$ contains a $U$-meager set of Loeb measure one. A cut is called bad if it is not good. Keisler and Leth discovered that most cuts are good and bad cuts are difficult to construct. In fact, they constructed bad cuts in some nonstandard universes under some extra set theoretic assumption beyond $ZFC$ such as $2^\omega < 2^{\omega_1}$. They proved also in [KL] that a bad cut must be a type two cut and a type two cut must have both uncountable cofinality and uncountable coinitiality. Given a cut $U$, the cofinality of $U$ is the cardinal

$$cof(U) = \min\{card(S) : S \subseteq U \land (\forall x \in U)\,(\exists y \in S)\,(x < y)\}$$

and the coinitiality of $U$ is the cardinal

$$coin(U) = \min\{card(S) : S \subseteq {}^*\mathbb{N} \smallsetminus U \land (\forall x \in {}^*\mathbb{N} \smallsetminus U)\,(\exists y \in S)\,(y < x)\},$$

where $card(S)$ denotes the external cardinality of $S$. The questions whether there exists a nonstandard universe in which there are no bad cuts or no type two cuts or no cuts $U$ with $cof(U) > \omega$ and $coin(U) > \omega$ are asked in [KL]. In [J1] the author showed that (1) bad cuts exist in some nonstandard universe (eliminating the need of the assumption $2^\omega < 2^{\omega_1}$), (2) in any $\omega_2$-saturated nonstandard universe there exist cuts $U$ with $cof(U) > \omega$ and $coin(U) > \omega$, (3) assuming $\mathfrak{b} > \omega_1$, i.e. every $B \subseteq \omega^\omega$ of cardinality $\leqslant \omega_1$ is eventually dominated by some $f \in \omega^\omega$, then every hyperfinite time line in any nonstandard universe has cuts $U$ with $cof(U) > \omega$ and $coin(U) > \omega$. Later Shelah [Sh] proved a surprising result that every hyperfinite time line in any nonstandard universe has cuts $U$ with $cof(U) = coin(U)$ without using any extra set theoretic assumption. Note that $cof(U) = coin(U)$ implies $cof(U) > \omega$ and $coin(U) > \omega$ by $\omega_1$-saturation. The main idea in the proof of the existence of type two cuts in §1 is the combination of Shelah's method of constructing cuts $U$ with $cof(U) > \omega$ and $coin(U) > \omega$ in [Sh] and Keisler-Leth's method of constructing type two cuts in [KL]. In the first half of the second section we prove that if the nonstandard universe satisfies the $\aleph_1$-isomorphism property, than there exist bad cuts in every hyperfinite time line. The proof uses a method from [JS]. In the second half of the second section we deal with very bad cuts (see definition below).

Suppose $U$ is a bad cut in some hyperfinite time line $\mathcal{H}$. By [KL, Proposition 4.5] $\mathcal{H}$ contains no $U$-meager set with positive Loeb measure. So if $S \subseteq \mathcal{H}$ is a $U$-meager set, then $S$ is either a non-Loeb measurable set or a Loeb measure zero set. A cut $U$ in $\mathcal{H}$ is called very bad if every $U$-meager set has Loeb measure zero. In the second section we prove that if the nonstandard universe satisfies the $\aleph_1$-isomorphism property, then for any cut $U$ except $U = H/\mathbb{N}$ in a hyperfinite time line $\mathcal{H}$, there exists a $U$-nowhere dense set $S \subseteq \mathcal{H}$ such that $S \not\subseteq A$ for any internal $A \subseteq \mathcal{H}$ with $\mu(A) \not\approx 1$, and $A \not\subseteq S$ for any internal $A \subseteq \mathcal{H}$ with $\mu(A) \not\approx 0$ (we then say that $S$ has outer Loeb measure one and inner Loeb measure zero). So if $U$ is a bad cut, then there is a non-Loeb measurable $U$-nowhere dense subset of $\mathcal{H}$. Hence $U$ is not very bad.

This paper is a sequel to [KL], [J1] and [Sh]. The reader is recommended to consult [CK] for background in model theory, to consult [CK], [L] or [SB] for background in nonstandard analysis, nonstandard universes and Loeb measure construction. In this paper we shall write $card(S)$ for the external cardinality of the set $S$ and write $|A|$ for the internal cardinality of $A$ when $A$ is an internal set. Let $^*\mathbb{R}$ denote the set of all real numbers in a given nonstandard universe $V$. For each $r \in {}^*\mathbb{R}$ we shall write $[r]$ for the greatest integer less than or equal to $r$. We call a number $r \in {}^*\mathbb{R}$ bounded if there is an $n \in \mathbb{N}$ such that $-n < r < n$. Otherwise we call $r$ unbounded (the use of the word *unbounded* here may avoid the confusion of using too much the word *infinite*). We call an $r \in {}^*\mathbb{R}$ infinitesimal if for any $n \in \mathbb{N}$ we have $-\frac{1}{n} < r < \frac{1}{n}$. We write $r \approx s$ if $r - s$ is an infinitesimal. We call a set infinite if it is externally infinite.

**Acknowledgements** The main results of this paper were obtained when the author was a visitor in University of Illinois at Urbana-Champaign during 94-95 year. The author very much appreciates the opportunity to work there in a friendly scholastic environment.

# 1 Type Two Cuts

Let's fix a nonstandard universe $V$ through out this section.

**Theorem 1** *There are type two cuts.*

**Proof:** The use of $^*\mathbb{R}$ in the following is not necessary and is only for convenience. One could replace any number $r$ in $^*\mathbb{R}$ by $[r] \in {}^*\mathbb{N}$ if he insists on working in a nonstandard model of Peano Arithmetic. In order to avoid multiple superscripts we write $\exp(a, b)$ for $a^b$ when $a, b \in {}^*\mathbb{R}$ and $a > 0$. First we construct sequences $\langle a_{n,\alpha} : \alpha < \lambda \rangle$ and $\langle b_{n,\alpha} : \alpha < \lambda \rangle$ for all $n \in \mathbb{N}$ simultaneously by a transfinite induction on ordinal $\lambda$ such that for any $n \in \mathbb{N}$ and $\alpha, \beta \in \lambda$ the following conditions are satisfied.

(a) $a_{n,\alpha}$ and $b_{n,\alpha}$ are positive and unbounded in ${}^*\mathbb{R}$.
(b) $a_{n,\alpha} < b_{n,\alpha}$.
(c) $\alpha < \beta \longrightarrow a_{n,\alpha} < a_{n,\beta}$.
(d) $\alpha < \beta \longrightarrow b_{n,\alpha} > b_{n,\beta}$.
(e) $a_{n,\alpha} = \exp(b_{n,\alpha}, 1/b_{n+1,\alpha}^3)$.
(f) $\alpha + 1 < \lambda \longrightarrow b_{n,\alpha+1} = \exp(b_{n,\alpha}, 1/b_{n+1,\alpha})$.
(g) $\alpha + 1 < \lambda \longrightarrow \exp(a_{n,\alpha}, b_{n+1,\alpha}) \leqslant a_{n,\alpha+1}$.

Suppose the construction is done up to stage $\lambda$. It is easy to see that for each $n \in \mathbb{N}$ the sequence $\langle a_{n,\alpha} : \alpha < \lambda \rangle$ is increasing, the sequence $\langle b_{n,\alpha} : \alpha < \lambda \rangle$ is decreasing and all $a_{n,\alpha}$'s are below all $b_{n,\alpha}$'s. For each $n \in \mathbb{N}$ let

$$J_{n,\lambda} = \bigcap_{\alpha < \lambda} \{x \in {}^*\mathbb{N} : a_{n,\alpha} < x < b_{n,\alpha}\}.$$

We shall show that if $J_{n,\lambda} \neq \emptyset$ for all $n \in \mathbb{N}$, then the inductive construction continues. So when the construction can not go further, there must be an $n \in \mathbb{N}$ such that $J_{n,\lambda} = \emptyset$. In this case, we shall use the sequence $\langle a_{n,\alpha} : \alpha < \lambda \rangle$ to define a type two cut.

Given any hyperinteger $H$, we choose a decreasing sequence $\langle d_n : n \in \mathbb{N} \rangle$ in ${}^*\mathbb{N} \smallsetminus \mathbb{N}$ such that $d_1 \leqslant H$ and $\exp(d_{n+1}, d_{n+1}^3) < d_n$ for each $n \in \mathbb{N}$. The sequence $\langle d_n : n \in \mathbb{N} \rangle$ exists by overspill principle. For the first step of the induction we choose $b_{n,0} = d_n$ and $a_{n,0} = \exp(d_n, 1/d_{n+1}^3)$. It is easy to see that for $\lambda = 1$ the conditions (c), (d), (f) and (g) are vacuously true and the conditions (b) and (e) are trivially true. For (a) since $d_{n+1} \leqslant \exp(d_n, 1/d_{n+1}^3)$, then $d_{n+1} \leqslant a_{n,0}$.

Suppose now the sequences $\langle a_{n,\alpha} : \alpha < \lambda \rangle$ and $\langle b_{n,\alpha} : \alpha < \lambda \rangle$ have been constructed such that for any $n \in \mathbb{N}$ and $\alpha, \beta < \lambda$ the conditions (a)—(g) are satisfied.

Case 1: $\lambda = \gamma + 1$ for some ordinal $\gamma$. For each $n \in \mathbb{N}$ let $b_{n,\lambda} = \exp(b_{n,\gamma}, 1/b_{n+1,\gamma})$ and let $a_{n,\lambda} = \exp(b_{n,\lambda}, 1/b_{n+1,\lambda}^3)$. We need to show that the sequences

$$\langle a_{n,\alpha} : \alpha < \lambda + 1 \rangle \text{ and } \langle b_{n,\alpha} : \alpha < \lambda + 1 \rangle$$

satisfy the conditions (a)—(g) with $\lambda$ replaced by $\lambda + 1$. Note that the conditions (b), (d), (e) and (f) are trivially true.

**Claim 1.1** The condition (g) is true.

Proof of Claim 1.1: First we have

$$\begin{aligned} b_{n+1,\gamma} &= \exp(b_{n+1,\lambda}, b_{n+2,\gamma}) \\ &> \exp(b_{n+1,\lambda}, 3) \\ &= b_{n+1,\lambda}^3 \end{aligned}$$

since $b_{n+1,\lambda} = \exp(b_{n+1,\gamma}, 1/b_{n+2,\gamma})$ and $b_{n+2,\gamma} > 3$. So we now have

$$\begin{aligned} a_{n,\lambda} &= \exp(b_{n,\lambda}, 1/b_{n+1,\lambda}^3) \\ &= \exp(\exp(b_{n,\gamma}, 1/b_{n+1,\gamma}), 1/b_{n+1,\lambda}^3) \\ &> \exp(\exp(b_{n,\gamma}, 1/b_{n+1,\gamma}), 1/b_{n+1,\gamma}) \\ &= \exp(\exp(b_{n,\gamma}, 1/b_{n+1,\gamma}^3), b_{n+1,\gamma}) \\ &= \exp(a_{n,\gamma}, b_{n+1,\gamma}) \end{aligned}$$

Hence the condition (g) is true.

It is easy to see that (c) follows from (g) and (a) follows from (b) and (c).

Case 2: $\lambda$ is a limit ordinal. If there exists an $n_0 \in \mathbb{N}$ such that $J_{n_0,\lambda} = \emptyset$, then stop and the construction is finished. Otherwise choose $c_n \in J_{n,\lambda}$ for each $n \in \mathbb{N}$. Let $b_{n,\lambda} = c_n$ and let $a_{n,\lambda} = \exp(b_{n,\lambda}, 1/b_{n+1,\lambda}^3)$. We need to check that the sequences

$$\langle a_{n,\alpha} : \alpha < \lambda + 1 \rangle \text{ and } \langle b_{n,\alpha} : \alpha < \lambda + 1 \rangle$$

satisfy the conditions (a)—(g) with $\lambda$ replaced by $\lambda + 1$. Note that (b), (d), (e), (f) and (g) are trivially true.

**Claim 1.2** The condition (c) is true.

Proof of Claim 1.2: Given any $\alpha < \lambda$. Since $c_{n+1} < b_{n+1,\beta}$ for any $\beta < \lambda$, we have

$$\exp(a_{n,\alpha}, c_{n+1}^3) < \exp(a_{n,\alpha}, (b_{n+1,\alpha} b_{n+1,\alpha+1} b_{n+1,\alpha+2})).$$

Now by (g) we have

$$\begin{aligned} &\exp(a_{n,\alpha}, (b_{n+1,\alpha} b_{n+1,\alpha+1} b_{n+1,\alpha+2})) \\ \leqslant\ & \exp(a_{n,\alpha+1}, (b_{n+1,\alpha+1} b_{n+1,\alpha+2})) \\ \leqslant\ & \exp(a_{n,\alpha+2}, b_{n+1,\alpha+2}) \\ \leqslant\ & a_{n,\alpha+3} < c_n. \end{aligned}$$

So $\exp(a_{n,\alpha}, c_{n+1}^3) < c_n$. Hence

$$a_{n,\alpha} < \exp(c_n, 1/c_{n+1}^3) = a_{n,\lambda}.$$

It is now obvious that (a) follows from (c). This ends the construction.

Suppose the construction halts at stage $\lambda$ for some ordinal $\lambda$. Then $\lambda$ muct be a limit ordinal and there exists an $n \in \mathbb{N}$ such that $J_{n,\lambda} = \emptyset$. We want to construct a type two cut $U$ from the sequences constructed above. Let

$$U = \{y \in {}^*\mathbb{N} : (\exists \alpha < \lambda)\,(y < \log(a_{n,\alpha}))\},$$

where log is the logarithmic function of base 2. Let

$$M = \{y \in {}^*\mathbb{N} : (\forall \alpha < \lambda)\,(y < b_{n+1,\alpha})\}.$$

**Claim 1.3** $\{y \in {}^*\mathbb{N} : (\forall \alpha < \lambda)\,(\log(a_{n,\alpha}) < y < \log(b_{n,\alpha}))\} = \emptyset$.

Proof of Claim 1.3: Suppose the claim is not true. Let $y \in {}^*\mathbb{N}$ such that

$$\log(a_{n,\alpha}) < y < \log(b_{n,\alpha})$$

for all $\alpha < \lambda$. Then for any $\alpha < \lambda$ we have

$$a_{n,\alpha} < 2^y < b_{n,\alpha}.$$

This contradicts that $J_{n,\lambda} = \emptyset$.

**Claim 1.4** $U$ is a cut.

Proof of Claim 1.4: It is easy to see that $\mathbb{N} \subseteq U$. We want to show that $U$ is closed under addition. For any $x \in U$ it suffices to show that $2x \in U$. Let $x < \log(a_{n,\alpha})$ for some $\alpha < \lambda$. Then

$$\begin{aligned}2x &< 2\log(a_{n,\alpha}) \\ &= \log(a_{n,\alpha})^2 \\ &< \log(\exp(a_{n,\alpha}, b_{n+1,\alpha})) \\ &\leqslant \log(a_{n,\alpha+1}).\end{aligned}$$

So $2x \in U$.

**Claim 1.5** $M(U) = M$.

Proof of Claim 1.5: Let $x \in M$. Given any $y \in U$, we want to show that $xy \in U$. Let $y < \log(a_{n,\alpha})$ for some $\alpha < \lambda$. Then

$$\begin{aligned} xy &< b_{n+1,\alpha} \log(a_{n,\alpha}) \\ &= \log(\exp(a_{n,\alpha}, b_{n+1,\alpha})) \\ &\leqslant \log(a_{n,\alpha+1}). \end{aligned}$$

So $xy \in U$. This shows that $M \subseteq M(U)$.

Let $x \in {}^*\mathbb{N} \smallsetminus M$. We want to find a $y \in U$ such that $xy \notin U$. By the definition of $M$ there is an $\alpha < \lambda$ such that $x > b_{n+1,\alpha}$. Let $y = [\log(a_{n,\alpha+1})] + 1$. Then $y \in U$. We now have

$$\begin{aligned} xy &> b_{n+1,\alpha} \log(a_{n,\alpha+1}) \\ &= (b_{n+1,\alpha}/b^3_{n+1,\alpha+1}) \log(\exp(a_{n,\alpha+1}, b^3_{n+1,\alpha+1})) \\ &= (b_{n+1,\alpha}/b^3_{n+1,\alpha+1}) \log(b_{n,\alpha+1}). \end{aligned}$$

Since

$$b^3_{n+1,\alpha+1} < \exp(b_{n+1,\alpha+1}, b_{n+2,\alpha}) = b_{n+1,\alpha},$$

we have $(b_{n+1,\alpha}/b^3_{n+1,\alpha+1}) > 1$. So $xy > \log(b_{n,\alpha+1})$. So $xy \notin U$. This shows that $M(U) \subseteq M$.

**Claim 1.6** $xM \neq U$ for any $x \in U$ and $x/M \neq U$ for any $x \in {}^*\mathbb{N} \smallsetminus U$.

Proof of Claim 1.6: Given any $x \in U$. We want to show that $xM \neq U$. Let $x < \log(a_{n,\alpha})$ for some $\alpha < \lambda$. For any $y \in M$ we have

$$\begin{aligned} xy &< b_{n+1,\alpha} \log(a_{n,\alpha}) \\ &= \log(\exp(a_{n,\alpha}, b_{n+1,\alpha})) \\ &\leqslant \log(a_{n,\alpha+1}) \end{aligned}$$

by the condition (g). So $xM \subseteq \{1, 2, \ldots, [\log(a_{n,\alpha+1})]\}$. Hence $xM \neq U$ because $[\log(a_{n,\alpha+1})] + 1 \in U \smallsetminus xM$.

Given any $x \in {}^*\mathbb{N} \smallsetminus U$. We want to show that $x/M \neq U$. By Claim 1.3 there is an $\alpha < \lambda$ such that $x > \log(b_{n,\alpha})$. For any $y \in M$ we have

$$\begin{aligned} x/y &> (\log(b_{n,\alpha}))/b_{n+1,\alpha} \\ &= \log(\exp(b_{n,\alpha}, 1/b_{n+1,\alpha})) \\ &= \log(b_{n,\alpha+1}). \end{aligned}$$

So $\{1, 2, \ldots, [\log(b_{n,\alpha+1})]\} \subseteq x/M$. Hence $x/M \neq U$ because $[\log(b_{n,\alpha+1})] \in x/M \smallsetminus U$.

Combining all those claims we have that $U$ is a type two cut. $\square$

**Remarks:** (1) In the definition of type one–type two cuts and in the proof of Theorem 1 we never use $\omega_1$-saturation. So type two cuts also exist in any non-$\omega_1$-saturated nonstandard universe or any nonstandard model of Peano Arithmetic. But for a countable nonstandard model of Peano Arithmetic the existence of type two cuts has a much easier proof.

(2) Since a cut $U$ with $cof(U) = coin(U)$ in a nonstandard universe may not be a type two cut, Theorem 1 is stronger than the result of Shelah in [Sh] mentioned in the introduction.

## 2 Bad Cuts and Very Bad Cuts

Let's recall the definitions. A cut $U$ in a hyperfinite time line $\mathcal{H} = \{1, 2, \ldots, H\}$ is called a good cut if $\mathcal{H}$ contains a $U$-meager set of Loeb measure one, where a set $X \subseteq \mathcal{H}$ is called $U$-meager if $X$ is the union of countably many nowhere dense sets under $U$-topology. $U$ is called bad if it is not good. A bad cut $U$ is called very bad if every $U$-meager set has Loeb measure zero. We show in this section that the $\aleph_1$-isomorphism property implies that there exist bad cuts and there are no very bad cuts. This means that for any nonstandard universe $V$ if $V$ satisfies the $\aleph_1$-isomorphism property, then there exist bad cuts and there are no very bad cuts in any hyperfinite time line $\mathcal{H}$ in $V$.

Let's introduce the $\kappa$-isomorphism property for any infinite regular cardinal $\kappa$. Given a nonstandard universe $V$. Let $\mathcal{L}$ be a first-order language. An $\mathcal{L}$-structure $\mathfrak{A} = (A; \ldots)$ is called internally presented (in $V$) if the base set $A$ is internal (in $V$) and the interpretation in $\mathfrak{A}$ of each predicate symbol or function symbol of $\mathcal{L}$ is internal (in $V$). Let's fix a nonstandard universe $V$. $V$ is said to satisfy the $\kappa$-isomorphism property if the following is true.

> For any first-order language $\mathcal{L}$ with $card(\mathcal{L}) < \kappa$ and for any two internally presented $\mathcal{L}$-structures $\mathfrak{A}$ and $\mathfrak{B}$, if $\mathfrak{A}$ and $\mathfrak{B}$ are elementarily equivalent, then $\mathfrak{A}$ and $\mathfrak{B}$ are isomorphic.

The $\kappa$-isomorphism property was suggested by Henson [H1]. Henson's form of the $\kappa$-isomorphism property is simple and elegant, but less applicable in some cases. In this section we use only an equivalent form of Henson's property in [JS, Main Theorem] stated in Lemma 2 below, in terms of the satisfiability of some second-order types. This equivalent form makes the use of the $\kappa$-isomorphism property easier in our case. See [H1], [H2], [J2] and [JK] for the existence of nonstandard universes satisfying the $\kappa$-isomorphism property.

**Lemma 2** *Let $\kappa$ be any infinite regular cardinal. Then the $\kappa$-isomorphism property is equivalent to the following:*

*For any first-order language $\mathcal{L}$ with $card(\mathcal{L}) < \kappa$, for any internally presented $\mathcal{L}$-structure $\mathfrak{A}$ and for any set of $\mathcal{L} \cup \{X\}$-sentences $\Gamma(X)$, where $X$ is a new n-ary predicate symbol not in $\mathcal{L}$, if $\Gamma(X) \cup Th(\mathfrak{A})$ is consistent, then $\Gamma(X)$ is satisfiable in $\mathfrak{A}$, i.e. there exists an n-ary relation $R \subseteq A^n$ where $A$ is the base set of $\mathfrak{A}$ such that $(\mathfrak{A}, R) \models \varphi(R)$ for every $\varphi(X) \in \Gamma(X)$.*

**Remark:** The original proof of [JS, Main Theorem] has a minor restriction on $\kappa$, e.g. $\kappa < \beth_\omega$. But this restriction can been easily removed by using $\kappa$-saturation. See [Sch].

We need also an equivalent form of the bad-ness of a cut from [KL]. An internal function $f$ with $dom(f) = \{1, 2, \ldots, L_f\}$ for some $L_f \in {}^*\mathbb{N}$ is called an internal sequence. Given a cut $U$ in $\mathcal{H}$, a strictly increasing internal sequence $f$ in $\mathcal{H}$ is called a crossing sequence of $U$ if for any $x \in U$ there exists a $y \in range(f) \cap U$ such that $x < y$. The following lemma is a part of [KL, Proposition 4.5].

**Lemma 3** *A cut $U$ is bad iff for any crossing sequence $f$ of $U$ the internal sum*

$$\sum_{m=1}^{L_f - 1} (f(m)/f(m+1))$$

*is unbounded.*

**Theorem 4** *The $\aleph_1$-isomorphism property implies that there exist bad cuts in every hyperfinite time line.*

**Proof:** Fix a nonstandard universe $V$ satisfying the $\aleph_1$-isomorphism property. Given any hyperfinite time line $\mathcal{H} = \{1, 2, \ldots, H\}$ in $V$, we want to show that there exist bad cuts in $\mathcal{H}$. First we define an internally presented structure $\mathfrak{A}$. Let

$$\mathcal{F} = \quad \{f : f \text{ is an increasing internal sequence from} \\ \{1, 2, \ldots, L_f\} \text{ for some } L_f \leqslant H \text{ to } \mathcal{H}\}.$$

Then $\mathcal{F}$ is internal. Define an internally presented structure

$$\mathfrak{A} = (\mathcal{H} \cup \mathcal{F} \cup {}^*\mathbb{R}; \mathcal{H}, \mathcal{F}, R, S, \leqslant, +, \cdot, n)_{n \in \mathbb{N}},$$

where $A = \mathcal{H} \cup \mathcal{F} \cup {}^*\mathbb{R}$ is the base set of $\mathfrak{A}$, $\mathcal{H}$ and $\mathcal{F}$ are unary relation, $R$ is a ternary relation such that $\langle a, b, f \rangle \in R$ iff $f \in \mathcal{F}$, $a \in dom(f)$ and $f(a) = b$, $S$ is a function from $\mathcal{F}$ to ${}^*\mathbb{R}$ such that for any $f \in \mathcal{F}$

$$S(f) = \sum_{m=1}^{L_f - 1} (f(m)/f(m+1)),$$

$\langle {}^*\mathbb{R}; +, \cdot, \leqslant \rangle$ is the real field in $V$, and $n$ is a constant of the structure for each $n \in \mathbb{N}$. Let $\mathcal{L}$ be the language of $\mathfrak{A}$. Note that the following first-order $\mathcal{L}$-sentences are true in $\mathfrak{A}$.

$$\theta_n = \exists x (\mathcal{H}(x) \wedge x \geqslant n \wedge \forall y (\mathcal{H}(y) \longrightarrow y \leqslant x))$$

for each $n \in \mathbb{N}$, and

$$\eta = \begin{array}{l} \forall f \forall x \forall y (\mathcal{F}(f) \wedge \mathcal{H}(x) \wedge \mathcal{H}(y) \wedge x < y \longrightarrow \\ \exists g (\mathcal{F}(g) \wedge range(g) = range(f) \cap [x, y])). \end{array}$$

where the formula $range(g) = range(f) \cap [x, y]$ is an abbreviation of the first-order $\mathcal{L}$-formula

$$\forall z (\exists u R(u, z, g) \leftrightarrow x \leqslant z \wedge z \leqslant y \wedge \exists u R(u, z, f)).$$

Let $X \notin \mathcal{L}$ be a unary predicate symbol. We define $\Gamma(X)$ to be the set of $\mathcal{L} \cup \{X\}$-sentences which contains exactly the following:

$$\varphi_1(X) = \forall x (X(x) \longrightarrow \mathcal{H}(x))$$

$$\varphi_2(X) = \forall x \forall y (x \leqslant y \wedge \mathcal{H}(x) \wedge X(y) \longrightarrow X(x))$$

$$\varphi_3(X) = \forall x \forall y (X(x) \wedge X(y) \longrightarrow X(x + y))$$

$$\psi_n = \forall f (\mathcal{F}(f) \wedge \forall x (X(x) \longrightarrow \exists y \exists z (R(y, z, f) \wedge X(z) \wedge x \leqslant z)) \longrightarrow S(f) \geqslant n)$$

for each $n \in \mathbb{N}$.

Note that the sentences $\varphi_1(X)$, $\varphi_2(X)$ and $\varphi_3(X)$ say that $X$ is a cut in $\mathcal{H}$. The sentences $\psi_n(X)$ for $n \in \mathbb{N}$ say that if $f$ is a crossing sequence of $X$, then the internal sum $S(f)$ is unbounded. So $\Gamma(X)$ describes that $X$ is a bad cut by Lemma 3. So if $\Gamma(X)$ is satisfiable in $\mathfrak{A}$, then $\mathcal{H}$ must contain a bad cut. By Lemma 2 it suffices to show that $\Gamma(X) \cup Th(\mathfrak{A})$ is consistent.

Let $\mathfrak{A}'$ be a countable elementary submodel of $\mathfrak{A}$. Then $Th(\mathfrak{A}') = Th(\mathfrak{A})$. If we can show that $\Gamma(X)$ is satisfiable in $\mathfrak{A}'$, then it is clear that $Th(\mathfrak{A}) \cup \Gamma(X)$ is consistent.

**Claim 4.1**  $\Gamma(X)$ is satisfiable in $\mathfrak{A}'$.

Proof of Claim 4.1: Let $A' = \mathcal{H}' \cup \mathcal{F}' \cup \mathbb{R}'$ be the base set of $\mathfrak{A}'$ and let $\mathcal{F}' = \{f_i : i \in \mathbb{N}\}$. We now inductively construct an increasing sequence $\langle a_i : i \in \mathbb{N} \rangle$ and a decreasing sequence $\langle b_i : i \in \mathbb{N} \rangle$ in $\mathcal{H}'$ such that for each $i \in \mathbb{N}$

(a) $a_i < b_i$,
(b) $2a_i < a_{i+1}$,
(c) $b_i/a_i$ is unbounded in $\mathbb{R}'$,
(d) If $f \in \mathcal{F}'$ such that

$$range(f) = range(f_i) \cap \{x \in \mathcal{H}' : a_i \leqslant x \leqslant b_i\},$$

if $S(f)$ is bounded in $\mathbb{R}'$ and if $L_f$ is unbounded, then there is a $k \in \{1, 2, \ldots, L_f\} \cap \mathcal{H}'$ such that $f(k) \leqslant a_{i+1}$ and $f(k+1) \geqslant b_{i+1}$ (or $f$ has a jump across the interval $(a_{i+1}, b_{i+1})$).

We show first that the claim follows from the construction. Let

$$U = \{x \in \mathcal{H}' : (\exists i \in \mathbb{N})\, (x \leqslant a_i)\}.$$

Then $\varphi_1(U)$ and $\varphi_2(U)$ are trivially true in $(\mathfrak{A}', U)$. The sentence $\varphi_3(U)$ is true in $(\mathfrak{A}', U)$ by the condition (b). Given any $f_i \in \mathcal{F}'$ such that $f_i$ is a crossing sequence of $U$. To show that $\psi_n(U)$ is true in $(\mathfrak{A}', U)$ for any $n \in \mathbb{N}$ we need only to show that $S(f_i)$ is unbounded. Suppose $S(f_i)$ is bounded. By the fact that $\eta$ is true in $\mathfrak{A}'$ there exists a $g \in \mathcal{F}'$ such that

$$range(g) = range(f_i) \cap \{x \in \mathcal{H}' : a_i \leqslant x \leqslant b_i\}.$$

Then $S(g)$ is also bounded because $S(g) \leqslant S(f_i)$. Since $f_i$ is a crossing sequence of $U$, $a_i \in U$ and $b_i \notin U$, then $g$ is also a crossing sequence of $U$. Hence $L_g$ is unbounded (since no finite sequence could be a crossing sequence of any cut). By the condition (d) we know that $g$ has a jump from $a_{i+1}$ to $b_{i+1}$, i.e. $g(k) \leqslant a_{i+1}$ and $g(k+1) \geqslant b_{i+1}$ for some $k \in dom(g)$. So $g$ can't be a crossing sequence of $U$, a contradiction.

We now do the inductive construction. Choose any $a_1$ and $b_1$ in $\mathcal{H}'$ such that $b_1/a_1$ is unbounded (for example, $a_1 = 1$ and $b_1 = H$). Suppose we have found $\langle a_i : i < k \rangle$ and $\langle b_i : i < k \rangle$ for some $k > 1$ such that they satisfy the conditions (a)—(d). We need to find $a_k$ and $b_k$. Let $g \in \mathcal{F}'$ be such that

$$range(g) = range(f_{k-1}) \cap \{x \in \mathcal{H}' : a_{k-1} \leqslant x \leqslant b_{k-1}\}.$$

Case 1: $S(g)$ is unbounded or $L_g$ is bounded. Simply let $a'_k = a_{k-1}$ and $b'_k = b_{k-1}$.

Case 2: $S(g)$ is bounded and $L_g$ is unbounded. Let $n \in \mathbb{N}$ be such that $S(g) < n$. Since $g$ is an element in $\mathfrak{A}'$ and $\mathfrak{A}' \preceq \mathfrak{A}$, then there is a $t$ in $\mathfrak{A}'$ such that

$$t = \min\{g(m)/g(m+1) : m \in \mathcal{H}' \wedge 1 \leqslant m < L_g\}.$$

Let $m_0 \in \mathcal{H}'$ and $m_0 < L_g$ be such that $t = g(m_0)/g(m_0 + 1)$. Then

$$t(L_g - 1) \leqslant \sum_{m=1}^{L_g - 1} (g(m)/g(m+1)) = S(g) \leqslant n.$$

So we have $g(m_0 + 1)/g(m_0) \geqslant (L_g - 1)/n$. Now let $a'_k = g(m_0)$ and $b'_k = g(m_0 + 1)$.

Clearly we have $b'_k/a'_k$ is unbounded. Let $a_k = 2a'_k$ and $b_k = b'_k - 1$. Then it is easy to see that $b_k/a_k$ is still unbounded. Now it is obvious that the sequences

$$\langle a_i : i < k+1 \rangle \text{ and } \langle b_i : i < k+1 \rangle$$

satisfy the conditions (a)—(d). □

**Remarks:** (1) We don't know if it is true that bad cuts exist in any nonstandard universe without the $\aleph_1$-isomorphism property.

(2) The $\aleph_1$-isomorphism property is equivalent to the $\aleph_0$-isomorphism property plus $\omega_1$-saturation (see [J3] and [Sch]). In fact every $n \in \mathbb{N}$ is definable in $\mathfrak{A}$. So it is only for convenience that we add constants $n$ into the structure $\mathfrak{A}$.

(3) Given any hyperinteger $L$ and $K$ in $\mathcal{H}$ such that $K/L$ is unbounded. Then we can make the bad cut $U$ sit between $L$ and $K$, i.e. $L \in U$ and $K \notin U$. To do this, just add $L$ and $K$ as constants of $\mathfrak{A}$, add the sentences $X(L)$ and $\neg X(K)$ to $\Gamma(X)$ and let $a_1 = L$, $b_1 = K$ at the beginning of the inductive construction. See [KL, Proposition 7.10] for the motivation of this remark.

Next we show that the $\aleph_1$-isomorphism property implies the non-existence of very bad cuts.

**Theorem 5** *The $\aleph_1$-isomorphism property implies that for any hyperfinite time line $\mathcal{H} = \{1, 2, \ldots, H\}$ and for any $c \in \mathcal{H}$ such that $c/H \approx 0$ there exists an $X \subseteq \mathcal{H}$ such that*

*(1) $X$ has outer Loeb measure one and inner Loeb measure zero,*

*(2) for any $x, y \in X$ if $x \neq y$, then $|x - y| \geqslant c$.*

**Proof:** We use same method as in the proof of Theorem 4. Let's fix a nonstandard universe $V$ satisfying the $\aleph_1$-isomorphism property. Let $\mathcal{P}$ be the set of all internal subsets of $\mathcal{H}$. So $\mathcal{P}$ is internal. Define an internally presented structure

$$\mathfrak{A} = (\mathcal{H} \cup \mathcal{P} \cup {}^*\mathbb{R}; \mathcal{H}, \mathcal{P}, \in, \mu, +, \cdot, \leqslant, c, n)_{n \in \mathbb{N}},$$

where $A = \mathcal{H} \cup \mathcal{P} \cup {}^*\mathbb{R}$ is the base set of $\mathfrak{A}$, $\mathcal{H}$ and $\mathcal{P}$ are unary relations, $\in$ is the natural membership relations between the elements of $\mathcal{H}$ and the elements of $\mathcal{P}$, $\mu$ is a function from $\mathcal{P}$ to ${}^*\mathbb{R}$ such that for any $A \in \mathcal{P}$, $\mu(A) = |A|/H$, $\langle {}^*\mathbb{R}, +, \cdot, \leqslant \rangle$ is the real field in $V$, $c$ and $n$ for each $n \in \mathbb{N}$ are constants. Let $\mathcal{L}$ be the language of $\mathfrak{A}$ and let $X$ be a new unary predicate not in $\mathcal{L}$. Let $\Gamma(X)$ be the set of $\mathcal{L} \cup \{X\}$-sentences which contains exactly the following:

$$\theta_1(X) = \forall x(X(x) \longrightarrow \mathcal{H}(x))$$

$$\theta_2(X) = \forall x \forall y(X(x) \wedge X(y) \wedge x \neq y \longrightarrow |x - y| \geqslant c)$$

$$\varphi_n(X) = \forall A(\mathcal{P}(A) \wedge X \subseteq A \longrightarrow \mu(A) > 1 - \frac{1}{n})$$

for each $n \in \mathbb{N}$ and

$$\psi_n(X) = \forall A(\mathcal{P}(A) \wedge A \subseteq X \longrightarrow \mu(A) < \frac{1}{n})$$

for each $n \in \mathbb{N}$. It is easy to see that $\theta_1(X)$ says that $X$ is a subset of $\mathcal{H}$, $\theta_2(X)$ says that any two different elements of $X$ have distance greater or equal to $c$, $\varphi_n(X)$ for all $n \in \mathbb{N}$ say that $X$ has outer Loeb measure one and $\psi_n(X)$ for all $n \in \mathbb{N}$ say that $X$ has inner Loeb measure zero. So we are done if we can show that $\Gamma(X)$ is satisfiable in $\mathfrak{A}$. By Lemma 2 we need only to show the consistency of $\Gamma(X) \cup Th(\mathfrak{A})$. Let $\mathfrak{A}'$ be a countable elementary submodel of $\mathfrak{A}$. It suffices to show that $\Gamma(X)$ is satisfiable in $\mathfrak{A}'$. Let $A' = \mathcal{H}' \cup \mathcal{P}' \cup \mathbb{R}'$ be the base set of $\mathfrak{A}'$ and let $\mathcal{P}' = \{A_n : n \in \mathbb{N}\}$. We want to construct sets $\{x_n \in \mathcal{H}' : n \in \mathbb{N}\}$ and $\{y_n \in \mathcal{H}' : n \in \mathbb{N}\}$ such that for each $n \in \mathbb{N}$

(a) $\mu(A_n) \ll 1 \longrightarrow x_n \notin A_n$,

(b) $\mu(A_n) \gg 0 \longrightarrow y_n \in A_n$,

(c) $x_n \notin (\bigcup_{m<n}\{x \in \mathcal{H}' : x_m - c \leqslant x \leqslant x_m + c\}) \cup \{y_m : m < n\}$ and $y_n \notin \bigcup_{m\leqslant n}\{x \in \mathcal{H}' : x_m - c \leqslant x \leqslant x_m + c\}$,

where $a \ll b$ means $a < b$ and $a \not\approx b$. Suppose we have found $\{x_m : m < n\}$ and $\{y_m : m < n\}$ such that (a), (b) and (c) are true up to stage $n-1$. Let

$$I_n = \bigcup_{m<n} \{x \in \mathcal{H}' : x_m - c \leqslant x \leqslant x_m + c\}.$$

Note that

$$\mu(I_n \cup \{y_m : m < n\}) \leqslant (2c+1)n/H + n/H \approx 0.$$

Thus

$$\mu(\mathcal{H}' \smallsetminus (I_n \cup \{y_m : m < n\})) \approx 1.$$

If $\mu(A_n) \approx 1$, then choose any $x_n \in \mathcal{H}' \smallsetminus (I_n \cup \{y_m : m < n\})$. If $\mu(A_n) \ll 1$, then

$$\mathcal{H}' \smallsetminus (I_n \cup \{y_m : m < n\} \cup A_n) \neq \emptyset.$$

So we can choose $x_n$ in above set. Now let

$$J_n = I_n \cup \{x \in \mathcal{H}' : x_n - c \leqslant x \leqslant x_n + c\}.$$

Then $\mu(J_n) \approx 0$. Let $y_n \in A_n \smallsetminus J_n$ if $\mu(A_n) \gg 0$ and $y_n \in \mathcal{H}' \smallsetminus J_n$ otherwise. This ends the construction.

Let $W = \{x_n : n \in \mathbb{N}\}$. It is clear that $(\mathfrak{A}', W) \models \theta_1(W)$. By the first part of (c) we have $(\mathfrak{A}', W) \models \theta_2(W)$. By (a) we have $(\mathfrak{A}', W) \models \varphi_n(W)$ for every $n \in \mathbb{N}$ and by (b) and the second part of (c) we have $(\mathfrak{A}', W) \models \psi_n(W)$ for every $n \in \mathbb{N}$. So $\Gamma(X)$ is satisfiable in $\mathfrak{A}'$. $\square$

**Corollary 6** *The $\aleph_1$-isomorphism property implies that there are no very bad cuts.*

**Proof:** Given any hyperfinite time line $\mathcal{H} = \{1, 2, \ldots, H\}$ and given a bad cut $U \subseteq \mathcal{H}$. Since the cut $H/\mathbb{N}$ is a good cut, then $U \neq H/\mathbb{N}$. Let $c \in H/\mathbb{N} \smallsetminus U$. It is easy to see that $c/H \approx 0$. Let $X \subseteq \mathcal{H}$ be the set obtained in Theorem 5. Then $X$ is a $U$-nowhere dense set because any interval of length less than $c$ contains at most one element of $X$. Obviously $X$ does not have Loeb measure zero because it has outer Loeb measure one. So $U$ is not very bad. $\square$

Department of Mathematics, College of Charleston
Charleston, SC 29424

Department of Mathematics, Rutgers University
New Brunswick, NJ 08903
*e-mail: rjin@math.rutgers.edu*